\documentclass{article}
\usepackage{amssymb,amsmath,amsthm,amssymb}
\usepackage{xy,xypic}
\input xy
\xyoption{all}
\title{Generalized Matric Massey Products for Graded Modules\footnote{Mathematics Subject Classification 14D15,13D10,13D02,13D07,14D22}}
\author{Arvid Siqveland\\Buskerud University College\\PO. Box 251\\N-3601 Kongsberg,Norway\\email:arvid.siqveland@hibu.no}

\define\Hom{\operatorname{Hom}}
\define\Ext{\operatorname{Ext}}
\define\ext{\operatorname{ext}}
\define\ker{\operatorname{ker}}
\define\Def{\operatorname{Def}}
\define\im{\operatorname{im}}
\define\gr{\operatorname{gr}}
\define\GrA{\operatorname{GradAlg}}
\define\ob{\operatorname{ob}}
\define\cl{\operatorname{cl}}
\define\C{\operatorname{C}}
\define\uR{\underline{R}}
\define\um{\underline{m}}
\define\un{\underline{n}}
\define\uu{\underline{u}}
\define\ut{\underline{t}}
\define\uv{\underline{v}}
\define\ux{\underline{x}}
\define\ualpha{\underline{\alpha}}
\define\Gr{{\operatorname{Gr}}}
\define\id{{\operatorname{id}}}
\define\ull{\underline{\ell}}
\define\Sch{\underline{\operatorname{Sch}}}
\define\Sets{\underline{\operatorname{Sets}}}
\define\Proj{\operatorname{Proj}}
\define\Spec{\operatorname{Spec}}
\define\Hilb{\operatorname{Hilb}}

\newtheorem{remark}{Remark}
\newtheorem{definition}{Definition}

\newtheorem{lemma}{Lemma}
\newtheorem{proposition}{Proposition}
\newtheorem{corollary}{Corollary}

\begin{document}
\maketitle
\begin{abstract}
The theory of generalized matric Massey products has been applied
for some time to $A$-modules $M$, $A$ a $k$-algebra. The main
application is to compute the local formal moduli $\hat{H}_M$,
isomorphic to the local ring of the moduli of $A$-modules. This
theory is also generalized to $\mathcal{O}_X$-modules $\mathcal{M}$,
$X$ a $k$- scheme. In these notes we consider the definition of
generalized Massey products and the relation algebra in any
obstruction situation (a differential graded $k$-algebra with
certain properties), and prove that this theory applies to the case
of graded $R$-modules, $R$ a graded $k$-algebra, $k$ algebraically
closed. When the relation algebra is algebraizable, that is the
relations are polynomials rather than power series, this gives a
combinatorial way to compute open (\'{e}tale) subsets of the moduli
of graded $R$-modules. This also gives a sufficient condition for
the corresponding point in the moduli of
$\mathcal{O}_{\Proj(R)}$-modules to be singular. The computations
are straight forward, algorithmic, and an example on the postulation
Hilbert scheme is given.
\end{abstract}

\section{Introduction}
The theory of generalized matric Massey products (GMMP) for
$A$-modules, $A$ a $k$-algebra, is given by Laudal in
\cite{SLN1183}, and applied to the theory of moduli of global and
local modules in \cite{Siqveland011},\cite{Siqveland012}. This
theory can obviously be applied also to the study of various Hilbert
Schemes, leading to GMMP for graded  $R$-modules $M$, $R$ a graded
$k$-algebra. \vskip0,2cm Let $(A^\bullet,d_\bullet)$ be a
differential graded $k$-algebra, and let
$\ualpha=\{\alpha_{e_1},\dots,\alpha_{e_d}\}$ be a set of elements
in $H^1(A^\bullet).$ For $\un\in(\mathbb{N}-\{0\})^d,$
$|\un|=n_1+\cdots+n_d=2$ we have the ordinary cup-products
$$\ualpha\otimes_k\ualpha\rightarrow h^2(A^\bullet)$$ given by
$<\ualpha;\un>=\underset{\underset{\um_i\in({\mathbb{N}\cup\{0\}})^d}{\um_1+\um_2=\un}}\sum\alpha_{m_1}\cdot\alpha_{m_2}$
For example
$$<\ualpha;(1,0,\dots,0,1)>=\alpha_{e_1}\cdot\alpha_{e_d}+\alpha_{e_d}\cdot\alpha_{e_1}
\text{ and }
<\ualpha;(1,1,,0\dots,0)>=\alpha_{e_1}\cdot\alpha_{e_2}.$$

By inductively adding elements $\alpha_{\um}\in A^1,$
$\um\in\overline{B}\subseteq(\mathbb{N}\cup\{0\})^d$ due to some
relations, we define the higher order generalized matric Massey
products $<\ualpha;\un>$, $\un\in
B^\prime\subseteq(\mathbb{N}\cup\{0\})^d,$ for some $\un$ of higher
order $|\un|$, provided $A^\bullet$ satisfies certain properties.
The inductive definition of GMMP is controlled at each step by the
relations between the monomials in an algebra $\hat{H}_{\ualpha}$
constructed in parallel. We call this algebra the relation algebra
of $\ualpha.$ It is interesting in it own right to study the GMMP
structure and the relation algebra of various sets of $\ualpha\in
h^1(A^\bullet).$ \vskip0,2cm

Deformation theory is introduced as a tool for studying local
properties of various moduli spaces. It is well known that the
prorepresenting hull of the deformation functor of a point $M$ in
moduli is the completion of the local ring in that point
\cite{Schlessinger68}. Consider a graded $R$-module $M$, $R$ a
graded $k$-algebra. Choose a minimal resolution $0\leftarrow
M\leftarrow L_\bullet$ of $M$ and consider the degree zero part
$\Hom_{R,0}^{\bullet}(L_\bullet,L_\bullet)$ of the Yoneda complex.
Then $(\Hom_{R,0}^{\bullet}(L_\bullet,L_\bullet),d_\bullet)$ is a
differential graded $k$-algebra. Let
$\ux^\ast=\{x_1^\ast,\dots,x_d^\ast\}\subseteq
H^1(\Hom_{R,0}^{\bullet}(L_\bullet,L_\bullet))\cong\Ext^1_{R,0}(M,M)$
be a $k$-basis. Then the relation algebra $\hat{H}_{\ux^\ast}$ is
isomorphic to the prorepresenting hull $\hat{H}_M$ of the (graded)
deformation functor $\Def_M,$ i.e.
$\hat{H}_{\ux^\ast}\cong\hat{H}_M.$ In addition to the definition of
the graded GMMP, this is the main result of the paper, implying that
general results about the GMMP gives local information about moduli.
In addition, we get the following result, telling us how GMMP on $R$
can be used to study the singular locus of sheaves on $\Proj(R)$:

\setcounter{proposition}{1}
\begin{proposition} Let $M=\Gamma_\ast(\mathcal M)$ for $\mathcal
M\in\C_R(\Spec(k))$. Then $H_{\mathcal{M}}$, the hull of
$\Def_{R,\mathcal{M}}$, is nonsingular if $H_M$, the hull of
$\Def_{R,M}$, is.
\end{proposition}

\setcounter{proposition}{0}

We conclude the paper with an explicit example.

\section{Classical graded theory}
\subsection{Notation}

We let $R=\underset {d\in\mathbb{Z}}\oplus R_d$ be a graded
$k$-algebra, $k$ algebraically closed of characteristic $0$, $R$
finitely generated in degree $1$. We let
$M=\underset{d\in\mathbb{Z}}\oplus M_d$ be a graded $R$-module, and
let $M(p)$ denote the twisted module of $M$ with grading
$M_d(p)=M_{p+d}.$

We will reserve the name $S$ for the free polynomial $k$-algebra
$S=k[x_1,\dots,x_n]$, so that $R$ is a quotient of some $S$ by some
homogeneous ideal $I$, that is $R=S/I$.

By $\ull$ we mean the category of local artinian  $k$-algebras with
residue field $k$. If $U\in\ob(\ull)$ we will use the notation
$\um_U$ for the maximal ideal in $U$. A surjective morphism
$\pi:U\twoheadrightarrow V$ in $\ull$ is called small if
$\ker\pi\cdot\um_U=0$. The ring of dual numbers is denoted
$k[\varepsilon]$ , that is
$k[\varepsilon]=k[\varepsilon]/(\varepsilon^2).$

If $V$ is a vector space, $V^\ast$ denotes its dual.

\subsection{Homomorphisms}
Classically, homomorphisms of graded $k$ algebras $R$ are
\textit{homogeneous of degree $0$}. This is also the case with morphisms of
graded $R$-modules. We might extend this definition by giving a grading to the
homomorphisms:$$\Hom_R(\underset{d\in\mathbb{Z}}\oplus M_d,\underset
{d\in\mathbb{Z}}\oplus
N_d)=\underset{d\in\mathbb{Z}}\oplus\Hom_{R,d}(M,N)$$ where
$\phi_d\in\Hom_{R,d}(M,N)\subseteq\Hom_R(M,N)$ has the additional
property $$\phi_d(M_p)\subseteq N_{p+d}=N_p(d).$$

\subsection{Construction of graded $S$-modules}

For the graded $R$-modules $M$ and $N$, $M\oplus N$ does certainly
not inherit a total grading by $(M\oplus
N)_d=\underset{{d^\prime+d''=d}}\oplus(M_{d'}\oplus N_{d''})$. Thus
the sentence "a free graded $R$-module" does simply not make any
sense. In this section we clarify how the grading is given.

Recall that if $M$, $N$ are graded $R$-modules, $f:N\rightarrow M$ a
homogeneous homomorphism (of degree $0$), then $\ker(f)$ and
$\im(f)$ are both graded submodules.

\begin{lemma} Let $N=\underset{d\in\mathbb{Z}}\oplus N_d$ be a graded $R$-module, $M$ any
$R$-module and \newline $f:M\twoheadrightarrow N$ a surjective $R$-module
homomorphism. Then
$$\gr(f^{-1})=\underset{d\in\mathbb{Z}}\oplus f^{-1}(N_d)$$ has a natural structure of
 \textbf{\textit{graded}} $R$-module.
\end{lemma}
\begin{remark} The proof of the above lemma is immediate, but it is
not always the case that $M\cong\gr(f^{-1})$ for some $f$. In fact,
this is equivalent with $M$ being graded.
\end{remark}

For the sake of simplicity, assume that
$M=\underset{d\in\mathbb{Z}}\oplus M_d$ is a finitely generated
graded module, generated by a finite number of homogeneous elements
$m_1,\dots,m_n$ of degrees $p_1,\dots p_n$ respectively. Then we
have a surjective homomorphism
$$\varepsilon:R(-p_1)\oplus R(-p_2)\oplus\cdots\oplus R(-p_n)\rightarrow M\rightarrow
0$$ sending $e_i$ of degree $p_i$ to $m_i$ (also of degree $p_i$).

We easily see that $R(-p_1)\oplus R(-p_2)\oplus\cdots\oplus
R(-p_n)\cong\gr(\varepsilon^{-1}),$ making $R(-p_1)\oplus
R(-p_2)\oplus\cdots\oplus R(-p_n)$ into a graded module. As the
kernel is also generated by a finite number of homogeneous elements,
say by $(g_1,\dots,g_l)$, where $g_i=\sum_j g_{ij}$ with $g_{ij}$
homogeneous in $R$, we have the following:

\begin{proposition} Every finitely generated, graded $R$-module $M$
has a minimal resolution on the form
$$\cdots\rightarrow\underset{i=1}{\overset{m_n}\oplus} R(-d_i^n)^{\beta_{n,i}}\rightarrow\cdots\rightarrow\underset{i=1}{\overset{m_1}\oplus}R(d_i^1)^{\beta_{1,i}}\rightarrow M\rightarrow 0$$
Conversely, given homogeneous elements $g_1,\dots,g_n\in R$ of
degrees $d_1,\dots,d_n$ respectively, then
$R(-d_1)\oplus\cdots\oplus R(-d_n)$ maps surjectively onto the
graded module $(g_1,\dots g_n)\subseteq R,$ and so is a graded
$R$-module.
\end{proposition}

\subsection{Families of graded modules}

A finitely generated graded $R$-module $M$ defines a coherent sheaf
$\tilde{M}$ of $\mathcal O_{\Proj(R)}$-modules. In the same way, an
ideal $I\subseteq R$ defines a sheaf of ideals on $\Proj(R)$, and
this gives a subscheme of $\Proj(R)$. Thus the study of various
moduli spaces is influenced by the study of graded $R$-modules.

Let us denote $\underline{R}=\Spec(R)$ for short, and let $X$ be a
scheme/$k$. Then a sheaf of graded $\mathcal O_{\uR\times_k
X}$-modules is an $\mathcal O_{\uR\times_k X}$-module $\mathcal G$
such that $\mathcal G(\uR\times_k U)$ is a graded
$R\otimes_k\mathcal O_X(U)$-module for every open $U\subseteq X.$ We
define the contravariant functor $\Gr_R:\Sch_k\rightarrow\Sets$ by

$$\Gr_R(X)=\{\text{Coherent graded }\mathcal O_{\uR\times_k X}\text{-modules }\mathcal G_X|\mathcal G_X\text{ is }X\text{-flat}\}/\cong.$$

The moduli spaces that we want to give results about, are the
schemes representing various restrictions of the functor
$\C_R:\Sch_k\rightarrow\Sets$ given by

$$\C_R(X)=\{\text{coherent }\mathcal O_{\Proj(R)\times_k X}\text{-modules }\mathcal F|\mathcal F\text{ is }X\text{ flat}\}/\cong. $$

The restrictions can be that $\mathcal F$ is an ideal sheaf with
fixed Hilbert polynomial $p(t)$. Then the above functor in the case
where $R=k[t_1,\dots,t_n]$, is the Hilbert functor
$\Hilb_{\mathbb{P}^n_k}^{p(t)}$. If $\mathcal F$ is locally free of
rank $r$ with fixed chern classes, we get $\mathcal
M(\mathbb{P}^n_k;c_1,\dots,c_r)$ and so forth.

The two functors $\Gr_R$ and $\C_R$ are usually not equivallent,
i.e. there exists graded modules $M$ such that
$\Gamma_\ast(\tilde{M})\ncong M$. However, $\Gamma_\ast(\mathcal
F)\cong\Gamma_\ast(\mathcal G)\Rightarrow\mathcal F\cong\mathcal G$,
and this, as we will see, is sufficient for given applications.

\vskip0,2cm

 In general, for a contravariant functor $\mathcal
F:\Sch\rightarrow\Sets$, and an element $x\in\mathcal F(\Spec(k))$,
we define the fiber functor $\mathcal
F_x:\Sch/k\times\{\underline{pt}\}\rightarrow\Sets$  from the
category of pointed schemes over $k$ to the category of sets  by
$$\mathcal F_x(X)=\{F\in\mathcal F(X)|\mathcal
F(\Spec(k)\overset{\underline{pt}}\rightarrow
X)(F)=F_{\underline{pt}}=x\}.$$ If $\mathcal F$ is represented by a
scheme $\mathbb M$, and if $x\in\mathbb M$ is a geometric point,
then the tangent space in this point is
$$(\um_x/\um_x^2)^{\ast}\cong\Hom_x(\Spec(k[\varepsilon],X)\cong\mathcal F_x(\Spec(k[\varepsilon])).$$
The fiber functors define covariant functors
$D_x:\ull\rightarrow\Sets$, $D_x(V)=\mathcal F_x(\Spec(V))$. In our
situation, we obtain the two deformation functors
$$D_M^{\Gr_R}=\Def_{R,M},\text{ }D_{\mathcal M}^{\C_R}=\Def_{R,\mathcal M}:\ull\rightarrow\Sets$$
given by

$$\Def_{R,M}(V)=\{\text{f.g. graded }R\otimes_k
V\text{-modules} M_V|M_V\text{ is }V\text{-flat},\text{
}M_{V,0}\cong M\}/\cong$$ and

$$Def_{R,\mathcal M}(V)=\{\text{coherent }\mathcal O_{\Proj(R\otimes_k V)}\text{-modules}\mathcal M_V|M_V\text{ is }V\text{flat},\text{ }\mathcal M_{V,0}\cong\mathcal M\}/\cong.$$
For the rest of this section, we assume that $\mathcal
{M}=\tilde{M}$. We will use the notations $\Def_{M}$ and
$Def_{\mathcal M}$ when no confusion is possible.

By definition, the tangent spaces of the postulated moduli spaces
are
$$\Def_M(k[\varepsilon])\cong\Ext_{R,0}^1(M,M)\text{ and } \Def_{\mathcal M}(k[\varepsilon])\cong\Ext_{\Proj(R)}^1(\mathcal{M},\mathcal{M})$$
respectively ($\Ext_{R,0}^1(M,M)$ will be defined below). For every
$V\in\ull$, the morphism\hskip0,2cm $\tilde{ }$\hskip0,2cm is
surjective with section $\Gamma_\ast$. That is, in the diagram
$$\xymatrix{\Def_{R,M}(V)\ar[r]_{\widetilde{ }}&\Def_{R,\mathcal M}(V)\ar@/_1pc/[l]_{\Gamma_{\ast}}},$$
$\Gamma_\ast(\mathcal F)^{\tilde{}}=\mathcal F.$

For a surjective small morphism $\pi:U\twoheadrightarrow V$ in
$\ull$, given a diagram
 \label{liftdiagram}
$$\xymatrix{M_U\in\Def_{M}(U)\ar[r]_{\widetilde{}}\ar[d]&\Def_{\mathcal
M}(U)\ni \tilde{M}_U
\ar@/_1pc/[l]_{\Gamma_\ast}\ar[d]\\\Gamma_\ast(\mathcal
F_V)\in\Def_{M}(V)\ar[r]^{\widetilde{}}&\Def_{\mathcal
M}(V)\ni\mathcal F_V\ar@/^1pc/[l]^{\Gamma_\ast},}$$

with $M_U$ mapping to $\Gamma_\ast(\mathcal F_V)$,
$\Gamma_\ast(\mathcal F_V)$ mapping to $\mathcal F_V$. Then it
follows by functoriality of $\tilde{}$ that $\tilde{M}_U$ is a
lifting of $\mathcal{F}_V$.

This has obvious consequences, and we will eventually  prove the
following:

\begin{proposition}\label{proposition2} Let $M=\Gamma_\ast(\mathcal M)$ for $\mathcal
M\in\C_R(\Spec(k))$. Then $H_{\mathcal{M}}$, the hull of
$\Def_{R,\mathcal{M}}$, is nonsingular if $H_M$, the hull of
$\Def_{R,M}$, is.
\end{proposition}

\section{Deformation theory}\label{DefTheorySect}
\subsection{Generalized Massey Products}\label{GMMPsect}
In this subsection, we consider a differential graded $k$-algebra
$(A^{\bullet},d_{\bullet})$ with certain properties. We will assume
that $0\in\mathbb{N}$, and for $\un\in\mathbb{N}^d$, we will use the
notation $|\un|=\sum_{i=1}^d n_i.$ For
$\ualpha=(\alpha_{e_1},\dots,\alpha_{e_d})\in (H^1(A^\bullet))^d$,
$d\in\mathbb{N}$, we will define some generalized Massey products
$<\ualpha;\um>\in H^2(A^\bullet)$, $\um\in B^{\prime}$, where
$B^{\prime}\subseteq \{\un\in\mathbb{N}^d:|\un|\geq 2\}$. Notice
that the Massey products may not be defined for all (if any)
$n\in\mathbb{N}^d.$ The idea is the following:

Let $\alpha_{e_1},\dots,\alpha_{e_d}$ be a set of $d$ elements in
$H^1(A^\bullet)$, let
$B_2^\prime=\{\un\in\mathbb{N}^d:|\un|=2\}$, and put $\overline{B}_1=\{\un\in\mathbb{N}^d:|\un|\leq 1\}$. The first order Massey
products are then the ordinary cup-products in $A^\bullet.$ That is
$$<\ualpha;\un>=\overline{y(\un)},\text{ }
y(\un)=\underset{\underset{|m_i|=1}{\um_1+\um_2=\un}}\sum\alpha_{\um_1}\cdot\alpha_{\um_2},
\text{ }\un\in B_2^\prime.$$

\begin{definition}

We will say that the Massey product is \textit{identically zero} if
\newline $y(\un)=0$.
\end{definition}

 The higher order Massey products are
defined inductively: Assume that the Massey products are defined for
$\un\in B_N^\prime$,
$B_N^\prime\subseteq\{\un\in\mathbb{N}^d:|\un|\leq N\},
N\in\mathbb{N}.$ For each $\um\in
B_{N}\subseteq B_N^\prime$, assume there
exists a fixed linear relation $l(\um)=\sum_{l=0}^{N-2}\sum_{\un\in
B_{2+l}^\prime}\beta_{\un,\um}<\ualpha;\un>=0,$ and choose an
$\alpha_{\um}\in A^1$ such that $d(\alpha_{\um})=l(\um).$ The set
$\{\alpha_{\um}\}_{\um\in \bar{B}_N}$, $\bar{B}_N=\bar{B}_{N-1}\cup
B_N$ is called \textbf{\textit{a defining system}} for the \textbf{\textit{Massey
products}}
$$<\ualpha;\un>=\overline{y(\un)},\text{ }
y(\un)=\sum_{|\um|\leq
N+1}\underset{\underset{\um_i\in\bar{B}_N}{\um_1+\um_2=\um}}\sum\beta_{\un,\um}^\prime\alpha_{\um_1}\cdot\alpha_{\um_2},\text{
}\um\in B_{N+1}^\prime,$$ where $\beta_{\un,\um}^\prime$ are chosen
linear coefficients for each pair $\um,\un$ such that
$<\ualpha;\un>\in h^2(A^{\bullet})$.

One way to construct Massey products, that is to construct the
relations $\beta$ given above, is the following:
\vskip0,2cm

Let $\alpha_{e_1},\dots,\alpha_{e_d}$ be a set of representatives
for $d$ elements in $h^1(A^\bullet).$ Let
$S_2=k[[u_1,\dots,u_d]]/\um^2=k[[\uu]]/\um^2,$ $R_3=k[[\uu]]/\um^3$,
$\bar{B}_1=\{\un\in\mathbb{N}^d:|n|\leq 1\},$
$B_2^\prime=\{\un\in\mathbb{N}^d:|n|=2\},$
$\bar{B}_2^\prime=\bar{B}_1\cup B_2^\prime.$

\begin{definition}
The first order Massey products are the ordinary
cup-products in $A^\bullet.$ That is, for $\un\in B_2^\prime:$
$$<\ualpha;\un>=\overline{y(\un)},\text{ }y(\un)=\underset{\underset{|m_i|=1}{\um_1+\um_2=\un}}\sum\alpha_{\um_1}\cdot\alpha_{\um_2}.$$
\end{definition}

Choose a $k$-basis $\{y_1^\ast,\dots,y_r^\ast\}$ for
$h^2(A^\bullet)$, and put $$f_j^2=\sum_{\un\in
B_2^\prime}y_j(<\ualpha;\un>)\uu^{\un},\text{ }j=1,\dots,r.$$ Let
$S_3=R_3/(f_1^2,\dots,f_r^2),$ and choose $B_2\subseteq B_2^\prime$
such that $\{\uu^{\un}\}_{\un\in B_2}$ is a monomial basis for
$$\um^2/\um^3+(f_1^2,\dots,f_r^2).$$ Put $\overline{B}_2=\overline{B}_1\cup
B_2.$ For each $\un\in\mathbb{N}^d$ with $|\un|\leq 3,$ we have a
unique relation in $S_3:$
$$\uu^{\un}=\sum_{\um\in\overline{B}_2}\beta_{\un,\um}\uu^{\um},$$
and for each $\um\in B_2,$ $b_{\um}=\underset{\un\in
B_2^\prime}\sum\beta_{\un,\um}<\ualpha;\un>=0.$ Choose for each
$\um\in B_2$ an $\alpha_{\um}\in A^1$ such that
$d(\alpha_{\um})=-b_{\um}.$ Put
$$R_4=k[[\uu]]/\um^4+\um(f_2,\dots,f_r^2).$$
Choose a monomial basis $\{\uu^{\un}\}_{\un\in B_3^\prime}$ for
$\um^3/\um^4+\um^3\cap(f_1^2,\dots,f_r^2)$ such that for each
$\un\in B_3^\prime,$ $\uu^{\un}=u_k\cdot\uu^{\um}$ for some $0\leq
k\leq d$ and some $\um\in B_2.$

\begin{definition}The set $\{\alpha_{\um}\}_{\um\in\overline{B}_2}$ is
called a defining system for the Massey products
$<\ualpha;\un>,\text{ }\un\in B_3^\prime.$ \end{definition}

Assume that the $k$-algebras

$$\xymatrix{R_{n+1}=k[[\uu]]/(\um^{n+1}+\um(f_1^{n-1},\dots,f_r^{n-1}))\ar[d]_{\pi_{n+1}^\prime}\\S_n=k[[\uu]]/(\um^n+(f_1^{n-1},\dots,f_r^{n-1}))}$$

and the sets $B_{n-1}$, $\overline{B}_{n-1}$, $B_{n}^\prime$,
$\{\alpha_{\um}\}_{\um\in\overline{B}_{n-1}}$ has been constructed
for $1\leq n\leq N$ according to the above, in particular

$$\begin{aligned}
&\ker\pi_{n+1}^\prime=\um^n+(f_1^{n-1},\dots,f_r^{n-1})/(\um^{n+1}+\um(f_1^{n-1},\dots,f_r^{n-1})\\
&\cong\um^n/(\um^{n+1}+\um^n\cap\um(f_1^{n-1},\dots,f_r^{n-1})\oplus(f_1^{n-1},\dots,f_r^{n-1})/\um(f_1^{n-1},\dots,f_r^{n-1}),
\end{aligned}$$
and we assume (by induction) that $\{\uu^{\un}\}_{\un\in
B_{n}^\prime}$ is a basis for
$$I_{n+1}=\um^n/(\um^{n+1}+\um^n\cap\um(f_1^{n-1},\dots,f_r^{n-1}).$$

Put $\overline{B}_{N+1}^\prime=\overline{B}_N\cup B_{N+1}^\prime$.
For each $\un\in B_{N+1}^\prime$ we have a unique relation in
$R_{N+1}:$
$$\uu^{\un}=\sum_{\um\in\overline{B}_{N+1}^\prime}\beta_{\um,\un}^\prime\uu^{\um}+\sum_j\beta_{\un,j}f_j^{n-1.}$$

\begin{definition}
The $N$`th order Massey products are
$$<\ualpha;\un>=\overline{y(\un)},\text{ }\un\in B_{N+1}^\prime,$$
$$y(\un)=\sum_{|\um|\leq N+1}\sum_{\underset{\um_i\in\overline{B}_N}{\um_1+\um_2=\um}}\beta_{\um,\un}^\prime\alpha_{\um_1}\cdot\alpha_{\um_2}.$$
\end{definition}

For this to be well defined, we need both that $y(\un)$ is a
coboundary, and that its class is independent of the choices of
$\alpha$. We will only consider algebras $A^\bullet$ that obey this,
and we call $A^\bullet$ \textit{\textbf{an obstruction situation algebra}}, or an
\textbf{\textit{OS-algebra}} for short.

 Put
$$\xymatrix{f_j^N=f_j^{N-1}+\underset{\un\in B_{N}^\prime}\sum
y_j(<\ualpha;\un>)\uu^{\un},\\
R_{N+2}=k[[\uu]]/(\um^{N+2}+\um(f_1^N,\dots,f_r^N)),\ar[d]_{\pi_{N+1}^\prime}\\
S_{N+1}=R_{N+1}/(f_1^N,\dots,f_r^N)=k[[\uu]]/(\um^{N+1}+(f_1^N,\dots,f_r^N))\overset{\pi_{N+1}}\rightarrow
S_N, }$$ pick a monomial basis $\{\uu^{\un}\}_{\un\in B_{N}}$ for
$\ker\pi_{N+1}$ such that $B_{N}\subseteq B_{N}^\prime,$ and put
$\overline{B}_{N}=\overline{B}_{N-1}\cup B_{N}.$ For each
$\un\in\mathbb{N}^d,$ $|\un|\leq N$ we have a unique relation in
$S_{N+1},$
$\uu^{\un}=\underset{\um\in\overline{B}_{N}}\sum\beta_{\un,\um}\uu^{\um},$
and for each $\un\in B_{N},$
$$b_{\un}=\sum_{l=0}^{N-1}\sum_{\un\in
B_{2+l}^\prime}\beta_{\un,\um}<\ualpha;\un>=0.$$

For each $\um\in B_{N}$, choose $\alpha_{\um}\in A^1$ such that
$d(\alpha_{\um})=-b_{\um}.$

\begin{definition} The set $\{\alpha_{\um}\}_{\um\in\overline{B}_N}$ is
called a defining system for the Massey products
$<\ualpha;\un>,\text{ }\un\in B_{N+1}^\prime.$
\end{definition}

Choose a monomial basis $\{\uu^{\un}\}_{\un\in B_{N+1}^\prime}$ for
$$\um^{N+1}/\um^{N+2}+\um^{N+2}+\um^{N+1}\cap\um(f_1^N,\dots,f_r^N)$$
such that for each $\un\in B_{N+1}^\prime,$
$\uu^{\un}=u_k\cdot\uu^{\um}$ for some $1\leq k\leq d$ and $\um\in
B_{N+1}.$ The construction then continues by induction.

\begin{definition} Let $(A^\bullet,d_\bullet)$ be a differential
graded OS $k$-algebra, and let
$\alpha_{1},\dots,\alpha_{d}=\underline{\alpha}$ be a set of
elements in $h^1(A^\bullet).$ Let $\{y_1^\ast,\dots,y_r^\ast\}$ be a
$k$-basis for $H^2(A^\bullet),$ and for $1\leq j\leq r$, let
$$f_j=\sum_{l=2}^\infty\sum_{\un\in B_l^\prime}y_j(<\ualpha;\un>)\uu^{\un}.$$
Then we define
$$\hat{H}_{\ualpha}=k[[u_1,\dots,u_d]]/(f_1,\dots,f_r)$$ and call it
\textbf{The relation algebra} of $\ualpha.$
\end{definition}

\subsection{Obstruction theory}
In this section, fix once and for all a minimal (graded) resolution
of the graded $R$-module M:

$$\xymatrix{0&\ar[l]M&L_0\ar[l]_{\varepsilon}&L_1\ar[l]_{\delta_1}&L_2\ar[l]_{\delta_2}&\cdots\ar[l]_{\delta_3}}$$

with $L_n\cong \underset{i=1}{\overset{m_n}\oplus}
R(-d_{i,n})^{\beta_{i,n}}.$ Consider a small surjective morphism
$\pi:U\twoheadrightarrow V$ in $\ull$, and let $M_V\in\Def_M(V).$
Then an element $M_U\in\Def_M(U)$ such that $\Def_M(\pi)(M_U)=M_V$
is called a lifting of $M_V$ to $U$.

\begin{lemma}
\label{liftlemma} To give a lifting $M_U\in\Def_M(U)$ of the graded
$R$-module $M$ is equivalent to give a lifting of complexes:

$$\xymatrix{0&M_U\ar[l]\ar[d]&L_0\otimes_k U\ar[d]\ar[l]_{\varepsilon^U}&L_1\otimes_k
U\ar[d]\ar[l]_{\delta_1^U}&L_2\otimes_k
U\ar[l]_{\delta_2^U}\ar[d]&\cdots\ar[l]_{\delta_3^U}\\
0&M\ar[l]&L_0\ar[l]^{\varepsilon}&L_1\ar[l]^{\delta_1}&L_2\ar[l]^{\delta_2}&\cdots\ar[l]^{\delta_3}}$$
We also have that $\varepsilon^U(l\otimes_k 1)\in M\otimes_k\um_U$,
$\delta_i^U(l\otimes_k 1)\in L_{i-1} \otimes_k\um_U$ for all $i\geq
1$, and that the top row is exact.
\end{lemma}

\begin{proof}
Because $U$ is artinian, its maximal ideal $\um_U$ is nilpotent. We
will prove the lemma by induction on $n$ such that $\um_U^n=0$, the
case $n=1$ being obvious. Assuming the result true for $n$, assume
$\um_U^{n+1}=0$ and put $V=U/\um_U^n$. Then the sequence of
$U$-modules $0\rightarrow I\rightarrow U\overset\pi\rightarrow
V\rightarrow 0$ with $I=\um_U^n=\ker\pi$ is exact with $\pi$ a small
morphism, and such that $\um_V^n=0$. Notice that $M_U\otimes_U V$ is
$V$-flat and that $(M_U\otimes_U V)\otimes_V k\cong M_U\otimes_U
k\cong M$ such that $M_V:=M_U\otimes_U V\in\Def_M(V).$ Also notice
that $M_U\otimes_U I\cong(M_U\otimes_U k)\otimes_k I\cong M\otimes_k
I.$ Thus tensorising $0\rightarrow I\rightarrow
U\overset\pi\rightarrow V\rightarrow 0$ over $U$ with $M_U$ we get
the exactness of the first vertical sequence in the diagram

$$\xymatrix{&0\ar[d]&0\ar[d]&0\ar[d]&0\ar[d]&\\
0&M\otimes_k I\ar[l]\ar[d]&L_0\otimes_k I\ar[l]_{\varepsilon\otimes\id}\ar[d]&L_1\otimes_k I\ar[l]_{\delta_1\otimes \id}\ar[d]&\L_2\otimes_k I\ar[l]_{\delta_2\otimes\id}\ar[d]&\cdots\ar[l]_{\delta_3\otimes\id}\\
0&M_U\ar[l]\ar[d]&L_0\otimes_k
U\ar[d]\ar[l]_{\varepsilon^U}&L_1\otimes_k
U\ar[d]\ar[l]_{\delta_1^U}&L_2\otimes_k
U\ar[l]_{\delta_2^U}\ar[d]&\cdots\ar[l]_{\delta_3^U}\\
0&M_V\ar[l]\ar[d]&L_0\otimes_k V\ar[l]^{\varepsilon^V}\ar[d]&L_1\otimes_k V\ar[l]^{\delta_1^V}\ar[d]&L_2\otimes_k V\ar[l]^{\delta_2^V}\ar[d]&\cdots\ar[l]^{\delta_3^V}\\
&0&0&0&0}$$

The exactness of the horizontal top row follows from exactness of
$I$ over $k$, the bottom row is exact by assumption, and the middle
row is constructed as follows: Let $\varepsilon^U$ be a lifting of
$\varepsilon^V$, which obviously exists. By assumption,
$\varepsilon^V(l\otimes 1)=\varepsilon(l)\otimes 1+h$, $h\in
M\otimes_k\um_V.$ Thus $\varepsilon^V(l\otimes
1)=\varepsilon(l)\otimes 1+h+v,$ $h\in\um_U,$ $v\in I,$ that is
$\varepsilon^V(l\otimes 1)=\varepsilon(l)\otimes 1+u,$ $u\in\um_U.$
The commutativity of the first rectangle then follows  from the fact
that $\um_U\cdot I=0,$ i.e. $\pi$ is a small morphism.

Now, choose a lifting $\tilde{\delta}_1^U$ of $\delta_1^V$. As
above, $\tilde{\delta}_1^U(l\otimes 1)\in L_0\otimes_k\um_U$ by the
induction hypothesis, and therefore commutes with
$\delta_1\otimes\id.$

For each generator $l$ of $L_1\otimes_k U$, choose an $x\in
L_0\otimes_k I$ such that
$(\varepsilon\otimes\id)(x)=\varepsilon^U(\tilde{\delta}_1(l)),$ and
put $\delta_1(l)=\tilde{\delta}_1(l)-x.$ Then
$$\varepsilon^U(\delta_1(l))=\varepsilon^U(\tilde{\delta}_1(l))-\varepsilon^U(x)=
\varepsilon^U(\tilde{\delta}_1(l))-(\varepsilon\otimes\id)(x)=0.$$
This gives the desired lifting, and we may continue this way with
$\delta_i^U$, $i>1$. We have proved that the middle sequence is a
complex.

Conversely, given a lifting of complexes as in the lemma, then
taking the tensor product over $U$  with $V=U/m_U^n$ in the top row,
we get a lifting as in the above diagram. By the induction
hypothesis, the bottom row is exact with $M_V=H_0(L.\otimes_k V)$ a
lifting of $M$.

Writing up the long exact sequence of the short exact sequence of
complexes, we have that the sequence in the middle is also exact,
$M_U=H_0(L.\otimes_k U)$ is flat over $U$ because $M_U\cong
M\otimes_k U$ as $k$-vector space implies that $M_U$ is $U$-free and
thus flat.
\end{proof}

As the category of graded $R$-modules is the (abelian) category of
representations of the graded $k$-algebra $R$, a homomorphism
$\phi:M\rightarrow N$ between the two graded $R$-modules $M$ and $N$
is by definition homogeneous of degree $0$. This implies that the
derived functors of $\Hom_R$ in the category of graded $R$-modules,
is the derived functors of $\Hom_{R,0}$, where $\Hom_{R,0}$ denotes
$R$-linear homomorphisms of degree $0$. Thus we use the notation
$\Ext^p_{R,0}(M,N)$.

We have fixed the minimal graded resolution $0\leftarrow M\leftarrow
L_{\bullet}$ of $M$, and we define the graded Yoneda complex by
$$(\Hom^{\bullet}_{R,0}(L_{\bullet},L_{\bullet}),\delta_{\bullet})$$
where $\Hom_{R,0}^p(L_{\bullet},L_{\bullet})=\underset{n\geq
p}\Pi\Hom_{R,0}(L_n,L_{n-p})$ and where the differential\newline
$\delta_p:\Hom_{R,0}^p(L_{\bullet},L_{\bullet})\rightarrow\Hom_{R,0}^{p+1}(L_{\bullet},L_{\bullet})$
is given by

$$\delta_p(\{\xi_n\}_{n\geq p})=\{\delta_n\circ\xi_{n-1} -(-1)^p\xi_n\circ\delta_{n-p}\}_{n\geq
p+1},$$

the composition given by $\xi\circ\delta(x)=\delta(\xi(x)).$ It is
 straight forward to prove:

\begin{lemma}
$H^n(\Hom^{\bullet}_{R,0}(L_{\bullet},L_{\bullet}))\cong\Ext_{R,0}^n(M,N),$
$n\geq 0.$
\end{lemma}

\begin{proposition}
Let $\pi:U\twoheadrightarrow V$ be a small morphism in $\ull$ with
kernel $I$. Let $M_V\in\Def_M(V)$ correspond to the lifting
$(L_{\bullet}\otimes_k V,\delta_{\bullet}^V)$ of the complex
$(L_{\bullet},\delta)$.  Then there is a uniquely defined
obstruction
$$o(M_V,\pi)\in\Ext^2_{R,0}(M,M)\otimes_k I$$ given in terms of the
$2$-cocycle
$o\in\Hom^{\bullet}_{R,0}(L_{\bullet},L_{\bullet})\otimes_k I,$ such
that $o(M_V,\pi)=0$ if and only if $M_V$ may be lifted to $U$.
Moreover, if $o(M_V,\pi)=0,$ then the set of liftings of $M_V$ to
$U$ is a principal homogeneous space over $\Ext^1_{R,0}(M,M).$
\end{proposition}

\begin{proof}
Because $L_i$ is free for each $i$, we can choose a lifting
$\tilde{\delta}^R_i$ making the following diagram commutative for
each $i$:
$$\xymatrix{L_{i-1}\otimes_k U\ar[d]&L_i\otimes_k U\ar[l]_{\tilde{\delta}^U_i}\ar[d]\\
L_{i-1}\otimes_k V&L_i\otimes_k V\ar[l]^{\delta^V_i}}.$$ As $\pi$ is
small, the composition
$\tilde{\delta}^U_i\circ\tilde{\delta}^U_{i-1}:L_i\otimes_k
U\rightarrow L_{i-2}\otimes_k U$ is induced by a unique morphism
$o_i:L_i\rightarrow L_{i-2}\otimes_k I,$ and so
$$o=\{o_i\}\in\Hom^2_{R,0}(L_{\bullet},L_{\bullet})\otimes_k I.$$
Also, $o$ is a cocycle, and
$$o(M_V,\pi)=\bar{o}\in\Ext_{R,0}^2(M,M).$$
Another choice $\tilde{\delta}^R$, leads to an
$o\in\Hom^2_{R,0}(L_{\bullet},L_{\bullet})\otimes_k I$ differing by
the image of  an element in
$\Hom^1_{R,0}(L_{\bullet},L_{\bullet})\otimes_k I$ such that
$o(M_V,\pi)$ is independent of the choice of liftings. This also
proves the only if part.

If $o=o(M_V,\pi)=0$, then there is an element
$$\xi\in\\Hom^1_{R,0}(L_{\bullet},L_{\bullet})\otimes_k I$$ such
that $o=-d_1(\xi).$ Put $\delta_i^U=\tilde{\delta}_i^U+\xi_i,$ and
one finds that $\delta_i^U\circ\delta_{i-1}^U=0.$ Thus it follows
from lemma \ref{liftlemma} that $M_V$ can be lifted to $U$.

For the last statement, given two liftings $M_U^1$ and $M_U^2$
corresponding to $(l_{\bullet}\otimes_k U, \delta_{\bullet}^{U,1})$
and $(l_{\bullet}\otimes_k U, \delta_{\bullet}^{U,2})$. Then their
difference induce morphisms
$\eta_i=\delta_i^{U,1}-\delta_i^{U,2}:L_i\rightarrow
L_{i-1}\otimes_k I$, and (for each choice of basis element in $I$)
$\eta=\{\eta_i\}\in\Hom^{\bullet}_{R,0}(L_{\bullet},L_{\bullet})$ is
a cocycle and thus defining $\bar\eta\in\Ext^1_{R,0}(M,M).$ This
gives the claimed surjection $$\{\text{ Liftings of }M_V\text{ to }
U\}\times\Ext^1_{R,0}(M,M)\twoheadrightarrow\{\text{ Liftings of
}M_V\text{ to } U\}.$$
\end{proof}
Notice that this proves that
$\Hom^\bullet_{R,0}(L_\bullet,L_\bullet)$ is an
OS-algebra.\vskip0,2cm
We now
combine the theory of Massey products and the theory of
obstructions. We let $\hat{H}_M=\hat{H}$ denote the prorepresenting
hull (the local formal moduli) of $\Def_M.$ All the way we will use
the notations and constructions in section \ref{GMMPsect}.

Pick a basis $$\{x_1,\dots,x_d\}\in\Ext^1_{R,0}(M,M)^\ast,$$  and a
basis $$\{y_1,\dots,y_r\}\in\Ext^2_{R,0}(M,M)^\ast.$$ Denote by
$\{x_i^\ast\}$ and $\{y_i^\ast\}$ the corresponding dual
bases.\vskip0,2cm

Put $S_2=k[u_1,\dots,u_d]/\um^2=k[[\uu]]/\um^2,$
$\overline{B}_1=\{\un\in\mathbb{N}^d:|n|\leq 1\}$. We set

$$\alpha_{\underline{0}}=\{d_i\},\text{ }\alpha_{e_j}=\{x_{j,i}^\ast\}.$$

 Let $t_{\hat{H}}$  and
$t_{\Def_E}$ denote the tangent spaces of $\hat{H}$ and $\Def_M$
respectively. A deformation $E_2\in\Def_E(S_2)$ corresponding to an
isomorphism
$$t_{\hat{H}}\rightarrow t_{\Def_E},$$ is represented by the lifting $$\{L_\bullet\otimes_k S_2,d_\bullet^{S_2}\}\text{ of }\{L_\bullet,d_\bullet\}$$
where $$d_\bullet^{S_2}|_{L_\bullet\otimes
1}=\sum_{\um\in\overline{B}_1}\alpha_{\um\bullet}\otimes\uu^{\um},$$
Now, put $\pi^\prime_3:R_3=k[[\uu]]/\um^3\rightarrow S_2,$ choose
$B_2^\prime$ as in section \ref{GMMPsect}  and put
$\overline{B}_2^\prime=\overline{B}_1\cup B_2^\prime.$ Then
$$o(E_2,\pi^\prime_3)=\cl\{d_i^{S_2}\circ d_{i-1}^{S_2}\}=\sum_{\un\in B_2^\prime}\overline{y(\un)}\otimes\uu^{\un}\in\Ext^2_{R,0}\otimes_k\ker(\pi^\prime_3),$$
with
$y(\un)=\underset{\underset{\um_i\in\overline{B}_1}{\um_1+\um_2=\un}}\sum\alpha_{\um_1,i}\circ\alpha_{\um_2,i-1}.$
This is to say $<\ux^\ast;\un>=\overline{y(\un)}$ for each $\un\in
B_2^\prime.$

Translating, we get
$$\begin{aligned}o(E_2;\pi_3^\prime)&=\sum_{\un\in B_2^\prime}<\ux^\ast;\un>\otimes_k\uu^{\un}
=\sum_{i=1}^r y_i\otimes(\sum_{\un\in
B_2^\prime}y_i(<\ux;\un>)\uu^{\un})\\&=\sum_{i=1}^r y_i\otimes
f_i^2.\end{aligned}$$ Following the construction in section
\ref{GMMPsect}, for each $\um\in B_2,$ we pick a $1$-cochain
$\alpha_{\um}\in\Hom_{R,0}^1(L_\bullet,L_\bullet)$ such that
$$d(\alpha_{\um})=-b_{\um}=-\sum_{\un\in B_2^\prime}\beta_{\un,\um}y(\un).$$
Then the family $\{\alpha_{\um}\}_{\um\in\overline{B}_2}$ is a
defining system for the Massey products
$$<\ux^\ast;\un>\text{, }\un\in B_3^\prime.$$
Define $d^{S_3}$ by $$d^{S_3}|_{L_i\otimes
1}=\sum_{\um\in\overline{B}_2}\alpha_{\um,i}\otimes\uu^{\um}.$$ Then
$(d^{S_3})^2=0,$ and so, by lemma \ref{liftlemma},
$\{L_\bullet\otimes_k S_3,d_i^{S_3}\}$ corresponds to a lifting
$E_3\in\Def_M(S_3).$ We continue by induction: Given a defining
system $\{\alpha_{\um}\}_{\um\in\overline{B}_N}$ for the Massey
products $<\ux^\ast;\un>$, $\un\in B_{N+1}^\prime,$ and assume
$d^{S_N}$ is defined by $d^{S_N}|_{L_i\otimes
1}=\sum_{\um\in\overline{B}_N}\alpha_{\um,i}\otimes\uu^{\um}.$ It
then follows that
$$\begin{aligned}
o(E_N,\pi_{N+1}^\prime)=\sum_{\un\in
B_{N+1}^\prime}(\sum_{\underset{\um_i\in\overline{B}_N}{\um_1+\um_2=\um}}\beta_{\un,\um}^\prime\alpha_{\um_1,i}\circ\alpha_{\um_2,i-1})\otimes\uu^{\un}
+\sum_{j=1}^r y_j^\ast\otimes f_j^N.\end{aligned}$$ For $1\leq j\leq
r$, letting $f_j^{N+1}=f_j^N+\underset{\un\in B_{N+1}^\prime}\sum
y_j(<\ux^\ast;\un>)$ as in section \ref{GMMPsect}, gives
$$o(E_N,\pi^\prime_{N+1})=\sum_j y_j^\ast\otimes f_j^{N+1}.$$
Dividing out by the obstructions, that is letting
$S_{N+1}=R_{N+1}/(f_1^{N+1},\dots,f_r^{N+1}),$ makes the obstruction
$0$, that is $\underset{\un\in
B_{N+1}^\prime}\sum\beta_{\un,\um}<\ux^\ast;\un>=0$ for each $\um\in
B_{N+1},$ such that the next order defining system can be chosen.

We have proved the following:

\begin{proposition}\label{proposition4}
Let $R$ be a graded $k$-algebra, $M$ a graded $R$-module. Let
$\{\ux^\ast\}=\{x_1^\ast,\dots,x_d^\ast\}\subseteq\Ext^1_{R,0}(M,M)=h^1(\Hom_{R,0}^\bullet(M,M))$
be a $k$-vector space basis. Then the relation algebra of
$\{\ux^\ast\}$ is isomorphic to the prorepresenting hull $\hat{H}_M$
of $\Def_M$, that is

$$\hat{H}_{\{\ux^\ast\}}\cong\hat{H}_M.$$
\end{proposition}

\begin{proof}
This follows directly from Schlessingers article
\cite{Schlessinger68}.
\end{proof}

\begin{proof}Of proposition \ref{proposition2}:
For each small morphism $\pi:U\twoheadrightarrow V,$ if $M_V$ is
unobstructed, so is $\tilde{M}_V$. Thus, if there are no relations
in $H_M$, there are none in $H_{\tilde{M}}$ either.
\end{proof}

\begin{proposition}
Let $R$ be a graded $k$-algebra, $[\mathcal{M}]$ a point in the
moduli space $\mathbb{M}$ of $\mathcal{O}_{\Proj(R)}$-modules corresponding to $\mathcal{M}$. If all
cup-products of $M=\Gamma_\ast(\mathcal{M})$ are identically zero,
then $\mathbb{M}$ is nonsingular in the point $[\mathcal{M}]$.
\end{proposition}

\begin{proof} We can choose all defining systems for $M$ equal to
zero so that there are no relations in $H_M$. The result then
follows from proposition \ref{proposition2}.
\end{proof}

\section{An example of an obstructed determinantal variety in the
postulation Hilbert scheme}

I would like to thank Jan Kleppe for introducing me to the theory of
postulation Hilbert schemes.

The \textit{\textbf{postulation Hilbert scheme}} is the scheme parameterizing graded $R$-algebras
with fixed Hilbert function.
The following example is given to me by
him. The theory is treated in \cite{KleppeRoig04}.\vskip0,2cm

Let $R=k[x_0,x_1,x_2,x_3]$, $k=\overline{k}$ and consider the two
$R$-matrices

$$G_I=\left(\begin{matrix}x_0&x_1&x_2&x_3^3\\x_2&x_0&x_1&x_2^3\end{matrix}\right)\text{, }
G_J=\left(\begin{matrix}x_0&x_1&x_2&x_3^3\\x_3&x_0&x_1&x_2^3\end{matrix}\right).$$

We let $I$ and $J$ be the ideals generated by the minors of $G_I$
and $G_J$ respectively. Because the the graded
modules $M_I=R/I$ and $M_J=R/J$ have equal betti-series, they belong to the same
component in $\GrA$. Thus if the dimension
of the tangent space of the two modules differ, the one with the
highest dimension necessarily has to be obstructed (meaning that it
correspond to a singular point).
Computing with Singular \cite{Singular}, we find
$\ext^1_{S,0}(M_I,M_I)=24,$ $\ext^1_{S,0}(M_J,M_J)=22.$ We then know
that the first is an example of an obstructed module.

Notice that a computer program (a library in Singular) can be made
for these computations. This will be clear in this example. However,
for large tangent space dimensions, it seems that the common
computers of today are too small.

In this section we will cut out the tangent space by a hyperplane
where the variety in question is obstructed. This will give readable
information about the relations in the point corresponding to the
variety, and the example will be possible to read.

We put

$$\begin{aligned}
s_1&=x_1^2-x_0x_2,\text{ }s_2=x_0x_1-x_2^2,\text{
}s_3=x_0^2-x_1x_2,\text{ }\\
s_4&=x_2^4-x_1x_3^3,\text{ }s_5=x_1x_2^3-x_0x_3^3,\text{
}s_6=x_0x_2^3-x_2x_3^3.\end{aligned}$$

Then $I=(s_1,\dots,s_6)$ and $M=M_I=R/I$ is given my the minimal
resolution

$$0\leftarrow M\leftarrow R\overset{d_0}\leftarrow R(-2)^3\oplus
R(-4)^4\overset{d_1}\leftarrow R(-2)^2\oplus
R(-5)^6\overset{d_2}\leftarrow R(-6)^3\leftarrow 0$$ with

$$d_0=\left(\begin{matrix}s_1&s_2&s_3&s_4&s_5&s_6\end{matrix}\right),$$
$$d_1=\left(\begin{matrix}x_0&-x_2&x_3^3&0&0&x_2^3&0&0\\-x_1&x_0&0&x_3^3&0&0&x_2^3&0\\
x_2&-x_1&0&0&x_3^3&0&0&x_2^3\\0&0&x_1&x_0&0&x_0&x_2&0\\
0&0&-x_2&0&x_0&-x_1&0&x_2\\0&0&0&-x_2&-x_1&0&-x_1&-x_0\end{matrix}\right),$$
$$d_2=\left(\begin{matrix}x_3^3&x_2^3&0\\0&-x_3^3&x_2^3\\-x_0&-x_2&0\\x_1&x_0&0\\
-x_2&-x_1&0\\0&-x_0&x_2\\0&x_1&-x_0\\0&-x_2&x_1\end{matrix}\right).$$
To compute a basis for $\Ext^1_{R,0}(M,M),$ we apply the functor
$\Hom_{R,0}(-,M)$, resulting in the sequence
$$M\overset{d_0^T}\rightarrow M(2)^3\oplus
M(4)^4\overset{d_1^T}\rightarrow M(2)^2\oplus
M(5)^6\overset{d_2^T}\rightarrow M(6)^3\rightarrow 0.$$ We then
notice that $d_0^T=0$ and so $\Ext^1_{R,0}(M,M)=(\ker d_1^T)_{0}.$

Programming in Singular \cite{Singular}, a basis for
$\Ext^1_{R,0}(M,M)$ is given by the columns in the following
$6\times 24$-matrix

$$\left(\begin{matrix}
0&0&0&0&0&0&0&0\\
0&0&0&0&0&0&0&0\\
0&0&0&0&0&0&0&0\\
x_0x_3^3&x_0x_2x_3^2&x_0x_2^2x_3&x_0x_2^3&x_0x_1x_3^2&x_0x_1x_2x_3&x_0x_1x_2^2&x_0x_1^2x_3\\
x_2x_3^3&x_2^2x_3^2&x_2^3x_3&x_2^4&x_1x_2x_3^2&x_1x_2^2x_3&x_1x_2^3&x_1^2x_2x_3\\
x_1x_3^3&x_1x_2x_3^2&x_1x_2^2x_3&x_1x_2^3&x_1^2x_3^2&x_1^2x_2x_3&x_1^2x_2^2&x_1^3x_3
\end{matrix}\right.$$

$$\begin{matrix}
0&-x_0x_1&-x_1^2&-x_1x_2&-x_1x_3&x_0x_2&x_1x_2&x_2^2\\
0&-x_0^2&-x_0x_1&-x_0x_2&-x_0x_3&0&0&0\\
0&0&0&0&0&-x_0^2&-x_0x_1&-x_0x_2\\
x_0^2x_3^2&x_0x_3^3&x_1x_3^3&x_2x_3^3&x_3^4&0&0&0\\
x_0x_2x_3^2&0&0&0&0&x_0x_3^3&x_1x_3^3&x_2x_3^3\\
x_0x_1x_3^3&0&0&0&0&0&0&0\end{matrix}$$

$$\left.\begin{matrix}
x_2x_3&0&-x_0x_3&x_1x_3&0&-x_2x_3&x_0x_3&x_1x_3\\
0&x_2x_3&-x_2x_3&0&x_1x_3&-x_1x_3&0&x_0x_3\\
-x_0x_3&x_1x_3&0&-x_2x_3&x_0x_3&0&-x_1x_3&0\\
0&0&x_2^3x_3&0&0&x_1x_2^2x_3&0&0\\
x_3^4&0&0&x_2^3x_3&0&0&x_1x_2^2x_3&x_2^3x_3\\
0&x_3^4&0&0&x_2^3x_3&0&0&x_1x_2^2x_3\end{matrix}\right).$$

Following the algorithm and notation given in section \ref{GMMPsect},
we compute cup-products. Of the 300 computed, 79 are identically
zero in the meaning that
$$0\equiv\alpha_{i_1}\circ\alpha_{i^\prime_2}+\alpha_{i^\prime_1}\circ\alpha_{i_2}:R(-3)^2\oplus
R(-5)^6\rightarrow R.$$ Of the remaining 221, 205 are zero in
cohomology, giving in total 16 nonzero cup-products. With respect to
a basis $\{\tilde{y}_i\}_{i=1}^{33}$ for $\Ext^2_{R,0}(M,M),$ these
products can be expressed by

$$\begin{matrix}
&v_{13}v_{23}=y_1,&v_{13}v_{24}=y_2\\
&v_{17}v_{22}=-y_1,&v_{17}v_{24}=y_3\\
&v_{18}v_{22}=-y_2,&v_{18}v_{23}=-y_3\\
&v_{19}v_{22}=-y_2,&v_{19}v_{23}=-y_3\\
&v_{20}v_{23}=y_1,&v_{20}v_{24}=y_2\\
&v_{21}v_{22}=-y_1,&v_{21}v_{24}=y_3\\
&v_{22}^2=y_1,&v_{22}v_{23}=y_2\\
&v_{22}v_{24}=-y_3,&v_{23}^2=y_3\\
\end{matrix}$$

Letting

$$\begin{aligned}
f_1&=v_{13}v_{23}-v_{17}v_{22}+v_{20}v_{23}-v_{21}v_{22}+v_{22}^2\\
f_2&=v_{13}v_{24}-v_{18}v_{22}-v_{19}v_{22}+v_{20}v_{24}+v_{22}v_{23}\\
f_3&=v_{17}v_{24}-v_{18}v_{23}-v_{19}v_{23}+v_{21}v_{24}-v_{22}v_{24}+v_{23}^2
\end{aligned}$$

we may conclude from proposition \ref{proposition4}:

\begin{proposition} The determinantal scheme given by the minors of
the matrix $G_I$ has first order relations given by its second order
local formal moduli $$\hat{H}/\um^3\cong
k[[v_1,\dots,v_{24}]]/((f_1,f_2,f_3)+\um^3).$$
\end{proposition}

From \cite{KleppeRoig04} it follows that the obstruction space for
$M=R/I$ is $H^2(M,M,R)$ where $M$ is considered as a graded
$R$-algebra. This $k$-vector space has dimension 3, and so we may
conclude:

\begin{corollary}
The determinantal scheme given by the minors of the matrix $G_I$ is
maximally obstructed.
\end{corollary}

We are now going to put most (21) of the variables above to zero.
That is, we choose the most interesting of the 24 variables above;
$t_1=v_{22}$, $t_2=v_{23}$, $t_3=v_{24}$, all others are put to
zero. We follow the algorithm given in section \ref{GMMPsect} and we
work in the Yoneda complex $\Hom^{\bullet}(L_\bullet,L_\bullet)$
where $L_\bullet$ denotes the $R$-free resolution of $M$ given
above.

Notice that for $\alpha=\{\alpha_i\}\in\Hom^1(L_\bullet,L_\bullet)$
it is always sufficient to have the the two leading morphisms
$\alpha_1:L_1\rightarrow L_0$, $\alpha_2:L_2\rightarrow L_1$. Also,
it is known that finding these by the methods below, they can always
be extended to the full complex, see \cite{Siqveland011}.

We find
$$\begin{aligned}\alpha_{e_1,1}&=\left(\begin{matrix}
-x_2x_3&-x_1x_3&0&x_1x_2^2x_3&0&0\end{matrix}\right),\\
\alpha_{e_1,2}&=\left(\begin{matrix}
-x_3&0&-x_2^2x_3&0&0&0&0&0\\
0&x_3&0&-x_2^2x_3&0&-x_2^2x_3&0&0\\
0&0&0&0&0&0&0&0\\
0&0&0&-x_3&0&0&0&0\\
0&0&0&0&0&0&0&0\\
0&0&-x_3&0&0&0&0&0\end{matrix}\right),\\
\alpha_{e_2,1}&=\left(\begin{matrix}
x_0x_3&0&-x_1x_3&0&x_1x_2^2x_3&0\end{matrix}\right),\\
\alpha_{e_2,2}&=\left(\begin{matrix}
0&-x_3&0&0&0&x_2^2x_3&0&0\\
0&0&0&0&-x_2^2x_3&0&0&0\\
-x_3&0&0&0&0&0&0&0\\
0&0&0&0&-x_3&0&0&0\\
0&0&x_3&0&0&0&0&0\\
0&0&0&0&0&0&0&0\end{matrix}\right),\\
\alpha_{e_3,1}&=\left(\begin{matrix}
x_1x_3&x_0x_3&0&0&x_2^3x_3&x_1x_2^2x_3\end{matrix}\right),\\
\alpha_{e_3,2}&=\left(\begin{matrix} 0&0&0&0&x_2^2x_3&0&x_2^2x_3&0\\
0&0&0&0&0&0&0&x_2^2x_3\\
0&-x_3&0&0&0&0&0&0\\
0&0&x_3&0&0&0&0&0\\
0&0&0&x_3&0&0&0&0\\
0&0&0&0&0&0&0&0
\end{matrix}\right).\end{aligned}$$

This given, it is an easy match to compute the cup products (or the
first order generalized Massey products). This is

$$\begin{aligned}<\uv^\ast;(2,0,0)>&=\alpha_{e_1,1}\cdot\alpha_{e_1,2}=x_3^2(x_2,-x_1,x_2^3,0,0,x_1x_2^2,0,0)\\
<\uv^\ast;(1,1,0)>&=\alpha_{e_1,1}\cdot\alpha_{e_2,2}
+\alpha_{e_2,1}\cdot\alpha_{e_1,2}=x_3^2(-x_0,x_2,-x_0x_2^2,0,0,-x_2^3,0,0)\\
<\uv^\ast;(1,0,1)>&=\alpha_{e_1,1}\cdot\alpha_{e_3,2}
+\alpha_{e_3,1}\cdot\alpha_{e_1,2}\\&=x_3^2(-x_1,x_0,-x_1x_2^2,-x_0x_2^2,-x_2^3,-x_0x_2^2,-x_2^3,-x_1x_2^2)\\
<\uv^\ast;(0,2,0)>&=\alpha_{e_2,1}\cdot\alpha_{e_2,2}=x_3^2(x_1,-x_0,x_1x_2^2,0,0,x_0x_2^2,0,0)\\
<\uv^\ast;(0,1,1)>&=\alpha_{e_2,1}\cdot\alpha_{e_3,2}
+\alpha_{e_3,1}\cdot\alpha_{e_2,2}=x_3^2(0,0,x_2^3,x_1x_2^2,0,x_1x_2^2,x_0x_2^2,0)\\
<\uv^\ast;(0,0,0)>&=\alpha_{e_3,1}\cdot\alpha_{e_3,2}=x_3^2(0,0,0,x_2^3,x_1x_2^2,0,x_1x_2^2,x_0x_2^2).\end{aligned}$$

As classes in cohomology, we find (as we already knew)

$$\begin{aligned}
<\uv^\ast;(2,0,0)>&=y_1,\text{ }<\uv^\ast;(1,1,0)>=y_2,\text{ }
<\uv^\ast;(1,0,1)>=-y_3,\\<\uv^\ast;(0,2,0)>&=y_3,\text{ }
<\uv^\ast;(0,1,1)>=0,\text{ } <\uv^\ast;(0,0,2)>=0.\end{aligned}$$

We put $$f_1^2=t_1^2,\text{ }f_2^2=t_1t_2,\text{
}f_3^2=t_2^2-t_1t_3,$$ and the \textit{restricted} local formal
moduli to the second order is

$$\hat{H}^\prime_M/\um^3=k[[t_1,t_2,t_3]]/(f^2_1,f_2^2,f_3^3).$$

We follow the algorithm given in section \ref{GMMPsect} further. We
choose a basis $B_2$ for $\um^2/(\um^3+(f_1^2,f_2^2,f_3^2))$, e.g.
$$B_2=\{(0,2,0),(0,1,1),(0,0,2)\},$$
and choose a third order defining system: Notice that by
$<\uv^\ast;\un>$ we mean a  \textit{representative} of the
cohomology class.
$$\begin{aligned}b_{(0,2,0)}&=<\uv^\ast;(0,2,0)>+<\uv^\ast;(1,0,1)>\\
&=x_3^2(0,0,0,-x_0x_2^2,-x_2^3,0,-x_2^3,-x_1x_2^2).\end{aligned}$$
And similarly,
$$b_{(0,1,1)}=x_3^2(0,0,x_2^3,x_1x_2^3,x_1x_2^2,0,x_1x_2^2,x_0x_2^2,0),$$
$$b_{(0,0,2)}=x_3^2(0,0,0,x_2^3,x_1x_2^2,0,x_1x_2^2,x_0x_2^2).$$
Writing up why these are cocycles, we find what $\alpha_{\un}$ to choose for
$d(\alpha_{\un})=-b_{\un}$:
$$\begin{aligned}\alpha_{(0,2,0),1}&=\left(\begin{matrix}0&0&0&0&0&-x_0x_2x_3^2\end{matrix}\right),\\
\alpha_{0,2,0),2}&=\left(\begin{matrix}0&0&0&0&0&0&0&0\\0&0&0&0&-x_2x_3^2&0&-x_2x_3^2&0\\
0&0&0&0&0&0&0&-x_2x_3^2\\
0&0&0&0&0&0&0&0\\0&0&0&0&0&0&0&0\\0&0&0&0&0&0&0&0\end{matrix}\right),\end{aligned}$$

$$\begin{aligned}\alpha_{(0,1,1),1}&=\left(\begin{matrix}0&0&0&-x_0x_2x_3^2&0&0\end{matrix}\right),\\
\alpha_{0,1,1),2}&=\left(\begin{matrix}0&0&0&0&0&0&0&0\\0&0&x_2x_3^2&0&0&0&0&0\\
0&0&0&x_2x_3^2&0&x_2x_3^2&0&0\\
0&0&0&0&0&0&0&0\\0&0&0&0&0&0&0&0\\0&0&0&0&0&0&0&0\end{matrix}\right),\end{aligned}$$

$$\alpha_{(0,0,2),1}=\left(\begin{matrix}0&0&0&0&0&x_2^2x_3^2\end{matrix}\right),\text{
}\alpha_{(0,0,2),2}\equiv 0.$$

Again, following the algorithm given in section \ref{GMMPsect}, we
choose a monomial basis $B_3^\prime=\{(0,2,1),(0,1,2),(0,0,3)\}$ for
$\um^3/\um^4+\um^3\cap\um(f_1^2,f_2^2,f_3^2),$ and we compute the
Massey products (notice again that the next to last expression is
the representative in the Yoneda complex of its cohomology class):

$$\begin{aligned}<\uv^\ast;(0,2,1)>&=\alpha_{(0,2,0)}\cup\alpha_{e_3}+\alpha_{(0,1,1)}
\cup\alpha_{e_2}+\alpha_{(0,0,2)}\cup\alpha_{e_1}\\
&=(0,0,-x_2^2x_3^3,-x_1x_2x_3^3,0,-x_1x_2x_3^3,-x_0x_2x_3^3,0)=0\end{aligned}$$

$$\begin{aligned}
<\uv^\ast;(0,1,2)>&=\alpha_{(0,1,1)}\cup\alpha_{e_3}+\alpha_{(0,0,2)}\cup\alpha_{e_2}\\
&\equiv 0\end{aligned}$$

$$<\uv^\ast;(0,0,3)>=\alpha_{(0,0,2)}\cup\alpha_{e_3}\equiv 0.$$

We now put $$f_1^3=f_1^2,\text{ }f_2^3=f_2^2,\text{ and
}f_3^3=f_3^2$$ and so $$\hat{H}^\prime/\um^4\cong
k[[t_1,t_2,t_3]]/(f_1^3,f_2^3,f_3^3)+\um^4.$$ We put
$B_3=B_3^\prime=\{(0,2,1),(0,1,2),(0,0,3)\},$ and the next order
defining system is easy to find, only one of the representations of
the elements $b_{\un}$ is different from zero:
$$b_{(0,2,1)}=x_0x_3^3(0,0,-x_1,-x_0,0,-x_0,-x_2,0)+(0,0,x_3^3s_2,x_3^3s_3,0,x_3^3s_3,0,0).$$
We choose
$$\begin{aligned}
\alpha_{(0,2,1),1}&=\left(\begin{matrix}0&0&0&x_0x_3^3&0&0\end{matrix}\right),\\
\alpha_{(0,2,1),2}&=\left(\begin{matrix}0&0&0&0&0&0&0&0\\
0&0&-x_3^3&0&0&0&0&0\\
0&0&0&-x_3^3&0&-x_3^3&0&0\\
0&0&0&0&0&0&0&0\\
0&0&0&0&0&0&0&0\\
0&0&0&0&0&0&0&0\end{matrix}\right).
\end{aligned}$$
The rest of the elements in the defining system is chosen
identically zero, and we put
$B_4^\prime=\{(0,2,2),(0,1,3),(0,0,4)\}$ and compute the fourth
order Massey products:
$$\begin{aligned}
<\uv^\ast;(0,2,2)>&=\alpha_{(0,1,2)}\cup\alpha_{e_2}+\alpha_{(0,2,1)}\cup\alpha_{e_3}\\
&+\alpha_{(0,2,0)}\cup\alpha_{(0,0,2)}+\alpha_{(0,1,1)}\cup\alpha_{(0,1,1)}\\
&+\alpha_{(0,0,3)}\cup\alpha_{e_1}\equiv 0,
\end{aligned}$$

$$\begin{aligned}
<\uv^\ast;(0,1,3)>&=\alpha_{(0,1,2)}\cup\alpha_{e_3}+\alpha_{(0,0,3)}\cup\alpha_{e_2}\\
&+\alpha_{(0,1,1)}\cup\alpha_{(0,0,2)}\equiv 0
\end{aligned}$$

$$<\uv^\ast;(0,0,4)>=\alpha_{(0,0,3)}\cup\alpha_{e_3}+\alpha_{(0,0,2)}\cup\alpha_{0,0,2}\equiv
0.$$

Now, put $f_i^4=f_i^3,$ $i=1,2,3.$ Because these then are
homogeneous of degree two, the next order defining systems involves
only fourth order Massey products, and these can all be chosen
identically zero. Then the fifth order Massey products involves
$\alpha_{\um_1}\cup\alpha_{\um_2}$ with at least one of $|\um_i|=3$.
We see that $\alpha_{(0,2,1)}\cup\alpha_{\um}\equiv 0$ for all $\um$
with $|\um|=2,$ and so all fifth order Massey products are zero.
Noting also that $\alpha_{(0,2,1)}\cup\alpha_{(0,2,1)}\equiv 0$, we
are ready to conclude:

\begin{proposition}
Let $f_1=t_1^2$, $f_2=t_1t_2$, $f_3=t_2^2-t_1t_3$. Then there exist
an open subset of the the component of $\GrA$, the moduli scheme of
graded $R$-algebras, containing the determinantal scheme
corresponding to the matrix $G_I$ such that its intersection with
the hyperplane $t_4=\cdots=t_{24}=0$ is isomorphic to
$$k[t_1,t_2,t_3]/(f_1,f_2,f_3)$$ with versal family
$$M_{(t_1,t_2,t_3)}=k[x_0,x_1,x_2,x_3]/I((t_1,t_2,t_3))$$ for
$\ut\in Z(f_1,f_2,f_3)$ with
$$\begin{aligned}I(t_1,t_2,t_3)&=(s_1-x_2x_3t_1+x_0x_3t_2+x_1x_3t_3,s_2-x_1x_3t_1+x_0x_3t_3,\\
&s_3-x_1x_3t_2,s_4+x_1x_2^2x_3t_1-x_0x_2x_3^2t_2t_3+x_0x_3^3t_2^2t_3,\\
&s_5+x_1x_2^2x_3t_2+x_2^3x_3t_3,s_6+x_1x_2^2x_3t_3-x_0x_2x_3^2t_2^2+x_2^2x_3^2t_3^2).\end{aligned}$$
\end{proposition}

\begin{remark}
When the local formal moduli with its formal family is algebraizable in this way,
we get an open subset of the moduli (at least \'{e}tale). Thus we get a lot more than just
the local formal information. The conditions for when $\hat{H}_M$ is algebraizable is an
interesting question.
\end{remark}

\bibliographystyle{plain}
\bibliography{alggeom05}
\end{document}